\documentclass{amsart}
\usepackage{amsmath,amssymb,graphicx}

\newtheorem{thm}{Theorem}[section]
\newtheorem{lem}[thm]{Lemma}
\newtheorem{prop}[thm]{Proposition}
\newtheorem{cor}[thm]{Corollary}

\theoremstyle{definition}
\newtheorem{exmp}[thm]{Example}

\begin{document}

\title{Noncrossing partitions in surprising locations}
\author[J.~McCammond]{Jon McCammond}
      \address{Dept. of Math.\\
               U. C. Santa Barbara\\
               Santa Barbara, CA 93106}
      \email{jon.mccammond@math.ucsb.edu}

\maketitle

\section{Introduction.}

Certain mathematical structures make a habit of reoccuring in the most
diverse list of settings.  Some obvious examples exhibiting this
intrusive type of behavior include the Fibonacci numbers, the Catalan
numbers, the quaternions, and the modular group.  In this article, the
focus is on a lesser known example: the noncrossing partition lattice.
The focus of the article is a gentle introduction to the lattice
itself in three of its many guises: as a way to encode parking
functions, as a key part of the foundations of noncommutative
probability, and as a building block for a contractible space acted on
by a braid group. Since this article is aimed primarily at
nonspecialists, each area is briefly introduced along the way.

The noncrossing partition lattice is a relative newcomer to the
mathematical world.  First defined and studied by Germain Kreweras in
1972 \cite{Kr72}, it caught the imagination of combinatorialists
beginning in the 1980s \cite{DeZa86}, \cite{Ed80}, \cite{Ed82},
\cite{EdSi94}, \cite{Ka02}, \cite{Pr83}, \cite{Si00}, \cite{SiUl91},
\cite{St97}, and has come to be regarded as one of the standard
objects in the field.  In recent years it has also played a role in
areas as diverse as low-dimensional topology and geometric group
theory \cite{Bessis}, \cite{Br01}, \cite{BrWa02}, \cite{Kr00},
\cite{Kr02} as well as the noncommutative version of probability
\cite{An00}, \cite{An02}, \cite{NiSp97}, \cite{Sp94}, \cite{Sp97},
\cite{Sp98}, \cite{Vo95}, \cite{Vo00}.  Due no doubt to its recent
vintage, it is less well-known to the mathematical community at large
than perhaps it deserves to be, but hopefully this short paper will
help to remedy this state of affairs.

\section{A motivating example.}
Before launching into a discussion of the noncrossing partition
lattice itself, we quickly consider a motivating example: the Catalan
numbers.  The Catalan numbers are a favorite pastime of many amateur
(and professional) mathematicians.  In addition, they also have a
connection with the noncrossing partition lattice
(Theorem~\ref{thm:nc-props}).

\begin{figure}
\includegraphics[scale=1.2]{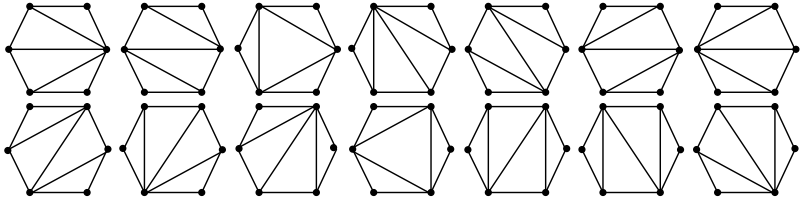}
\caption{The fourteen triangulations of a hexagon.\label{fig:triangulations}}
\end{figure}

\begin{figure}
\includegraphics[width=4in]{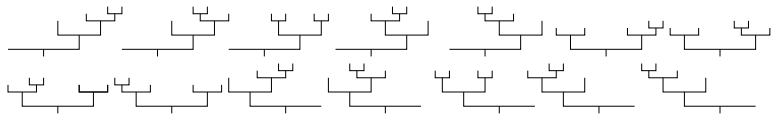}
\caption{The fourteen rooted trivalent plane trees with six leaves and
  ten vertices.\label{fig:trees}}
\end{figure}

\begin{exmp}[Catalan numbers]
The Catalan numbers are the numbers $C_n$ given by
\[C_n = \frac{1}{n+1}\left(\begin{array}{c}2n\\n\end{array}\right),\] 
\noindent 
and they have a number of different interpretations.  See, for
example, Richard Stanley's list of more than one hundred distinct ways
in which this sequence arises \cite{catalan}.  Some of the most common
interpretations are as the number of triangulations of an
$(n+2)$-sided polygon (illustrated in
Figure~\ref{fig:triangulations}), as the number of binary
parenthesizations of a string of $n+1$ letters
(Figure~\ref{fig:associate}), or as the number of rooted trivalent
plane trees with $2n+2$ vertices (Figure~\ref{fig:trees}).  With these
examples to whet the reader's appetite, we direct the interested
reader to \cite{catalan} and \cite[Exercise 6.19]{EC2}, and we
continue on to a description of noncrossing partitions.
\end{exmp}

\begin{figure}
$\begin{array}{cc}
(1 (2 (3 (4 5)))) &  (1 (2 ((3 4) 5))) \\
(1 ((2 3) (4 5))) &  (1 ((2 (3 4)) 5))\\
(1 (((2 3) 4) 5)) &  ((1 2) (3 (4 5)))\\
((1 2) ((3 4) 5)) &  ((1 (2 3)) (4 5))\\
((1 (2 (3 4))) 5) &  ((1 ((2 3) 4)) 5)\\
(((1 2) 3) (4 5)) &  (((1 2) (3 4)) 5)\\
(((1 (2 3)) 4) 5) &  ((((1 2) 3) 4) 5)\\
\end{array}$
\caption{The fourteen ways to asssociate five numbers.\label{fig:associate}}
\end{figure}

\section{Noncrossing partitions.}
We are now ready to define a noncrossing partition.  Following
traditional combinatorial practice we use $[n]$ to denote the set
$\{1,\ldots,n\}$.

\begin{figure}[ht]
\includegraphics[scale=1.5]{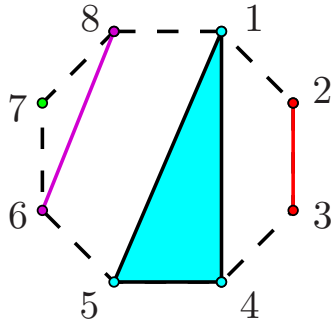}
\caption{A noncrossing partition of the set $[8]$.\label{fig:ncross}}
\end{figure}

\medskip
\noindent
{\bf Noncrossing partitions.}  Recall that a \emph{partition} of a set
is a collection of pairwise disjoint subsets whose union is the entire
set and that the subsets in the collection are called \emph{blocks}.
A \emph{noncrossing partition} $\sigma$ is a partition of the vertices
of a regular $n$-gon (labeled by the set $[n]$) so that the convex
hulls of its blocks are pairwise disjoint.  Figure~\ref{fig:ncross}
illustrates the noncrossing partition
$\{\{1,4,5\},\{2,3\},\{6,8\},\{7\}\}$.  The partition
$\{\{1,4,6\},\{2,3\},\{5,8\},\{7\}\}$ would be crossing.

\begin{figure}[ht]
\includegraphics[scale=1.3]{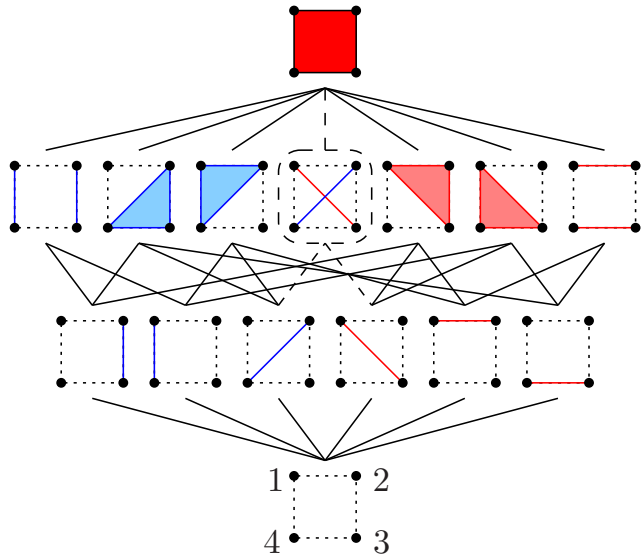}
\caption{The figure shows the partition lattice for $n=4$. If the
  vertex surrounded by a dashed line is removed, the result is the
  noncrossing partition lattice for $n=4$.\label{fig:nc4}}
\end{figure}

Given partitions $\sigma$ and $\tau$ of $[n]$ we say that $\sigma <
\tau$ if each block of $\sigma$ is contained in a block of $\tau$.
This ordering on the set of all partitions of $[n]$ defines a
partially ordered set called the \emph{partition lattice} and is
usually denoted $\Pi_n$.  When restricted to the set of noncrossing
partitions on $[n]$, it called the \emph{noncrossing partition
  lattice} and denoted $NC_n$.  The poset $\Pi_4$ is shown in
Figure~\ref{fig:nc4}.  For $n=4$, the only difference between the two
posets is the partition $\{\{1,3\},\{2,4\}\}$, which is not
noncrossing.

\medskip

The poset of noncrossing partitions has a number of nice combinatorial
properties that we record in the following theorem:

\begin{thm}\label{thm:nc-props}
For each $n$ the poset $NC_n$ is a graded, bounded lattice with
Catalan many elements ($C_n$ to be exact).  In addition, it is
self-dual and locally self-dual.
\end{thm}

\noindent
The fact that the number of noncrossing partitions is a Catalan number
is part of the lore of the Catalan numbers, and we refer the reader
again to \cite{catalan} and \cite{EC2}.  For the other properties we
now review their definitions.

A poset is \emph{bounded} if it has both a minimum element and a
maximum element.  For the noncrossing partition lattice, the discrete
partition (i.e., the one in which each block contains a single
element) and the partition with a single block fulfill these roles.  A
\emph{chain} in a poset is a subset in which any two elements are
comparable, its \emph{length} is one less than the size of this
subset, and a \emph{maximal chain} is a chain that is not properly
contained in any larger chain.  A poset in which any two maximal
chains have the same length is called \emph{graded}.  In any graded
bounded poset there is a height function that keeps track of the level
in which elements are contained.  The level of an element $x$ can be
defined as the number of elements strictly below $x$ in any maximal
chain containing $x$.  In $\Pi_n$ or $NC_n$, for example, the height
of a partition is $n$ minus the number of blocks.  For later use, we
also note that if $\sigma < \tau$ in $\Pi_n$ or $NC_n$ and if their
heights differ by one, then there are blocks $B$ and $B'$ in $\sigma$
whose union is a block in $\tau$.  Moreover, all other blocks in
$\sigma$ and $\tau$ are identical.  In this situation we say that
$\tau$ \emph{covers} $\sigma$.  A chain in which each element covers
the previous one is called a \emph{saturated chain}.

A poset is a \emph{lattice} if each pair of elements has a least upper
bound and a greatest lower bound.  The greatest lower bound $\sigma
\wedge \tau$ of partitions $\sigma$ and $\tau$ of $[n]$ is simply the
largest refinement of the two partitions.  In other words, define $i$
and $j$ to be in the same block of $\sigma \wedge \tau$ if and only if
they lie in the same blocks in both $\sigma$ and $\tau$.  It is now
easy to see that this is a lower bound for $\sigma$ and $\tau$ and
that this is greater than any other lower bound.  To find the least
upper bound $\sigma \vee \tau$ of $\sigma$ and $\tau$ from $NC_n$,
superimpose the convex hulls of all the blocks for $\sigma$ and for
$\tau$ and then take the convex hulls of the connected components that
result.

Finally, a poset is \emph{self-dual} if there is an order-reversing
bijection from it to itself, and it is \emph{locally self-dual} if
this is true for each of its intervals.  Recall that an interval
$[x,y]$ in a poset is simply the subposet containing all the elements
greater than or equal to $x$ and less than or equal to $y$.  Since
this property is easier to establish once we make the connection to
the symmetric groups, we postpone its proof until the next section.

\section{Symmetric groups.}
Before connecting the noncrossing partition lattice with
low-dimensional topology and noncommutative probability, it will be
helpful to establish first its close connection with the symmetric
group.

\medskip
\noindent
{\bf Symmetric groups.}  A \emph{permutation} of a set $X$ is a
bijection from $X$ to itself.  We use $S_n$ to signify the group of
all permutations of the set $[n]$ under function composition. (We
refer to $S_n$ as the \emph{symmetric group} on $n$ elements.)  There
are two natural generating sets for a group of permutations.  If the
underlying set is unordered, then the most natural generating set is
the set of permutations that interchange two elements and leave the
rest fixed (ususally called \emph{transpositions} or
\emph{two-cycles}).  If, on the other hand, the underlying set itself
has a natural linear ordering, as is the case for $[n]$, then the
smaller generating set using only adjacent transpositions $(i,i+1)$ is
often preferred.

One of the more mysterious properties of permutations for students in
an abstract algebra class is the fact that they are naturally
classified as either even or odd.  Since many of the standard proofs
of this are unenlightening and since there is an elementary geometric
proof, we digress slightly to present it.  The following
result is, of course, the heart of the matter:

\begin{thm}\label{thm:sym}
Every product of transpositions that equals the identity permutation
has an even number of factors.
\end{thm}

\begin{figure}[ht]
\begin{tabular}{ccc}
\begin{tabular}{c}
\includegraphics[scale=1]{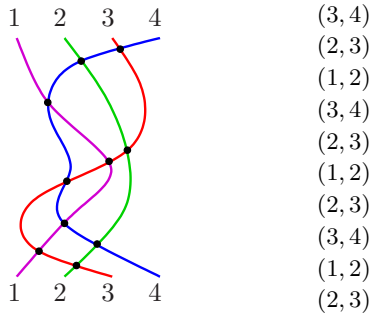}
\end{tabular}
& \hspace*{1cm} &
\begin{tabular}{c}
\small$(3,4)$\\
\small$(2,3)$\\
\small$(1,2)$\\
\small$(3,4)$\\
\small$(2,3)$\\
\small$(1,2)$\\
\small$(2,3)$\\
\small$(3,4)$\\
\small$(1,2)$\\
\small$(2,3)$\\
\end{tabular}
\end{tabular}
\caption{Illustration of the proof of Theorem~\ref{thm:sym}.\label{fig:sym}}
\end{figure}

\begin{proof}
We first prove the theorem in the case where all of the transpositions
involved are adjacent ones.  Given a sequence of adjacent
transpositions, we can draw a set of curves as illustrated in
Figure~\ref{fig:sym}.  The intersections correspond to the
transpositions as follows.  Arrange the transpositions from top to
bottom in the order they are applied.  In the example, the first
transposition is $(3,4)$ and the curves on the left are drawn so that
highest crossing occurs between the strands that are third and fourth
from the left.  The next transposition is $(2,3)$ so the curves are
drawn so that the next highest crossing occurs between the strands
that are currently second and third from the left.  This procedure can
be used to convert any sequence of adjacent transpositions into smooth
descending curves.  (This procedure can also be reversed for smooth
descending curves in general position.)

If the product is indeed the identity permutation, then the curve that
starts at position $i$ also ends at position $i$.  Call this the $i$th
curve.  Now we simply change perspectives.  Each intersection is an
intersection between two curves, say the $i$th and $j$th curves.
Since the $i$th and $j$th curves originally occur in a particular
order and they also return to the same order, these two curves must
intersect an even number of times.  Adding up the intersections
according to the curves involved gives the result.  To convert this to
a result about products of arbitrary transpositions, it remains only
to note that every transposition can be written as a product of an odd
number of adjacent transpositions, so products over the larger
generating set can be converted to products over the smaller
generating set without changing the parity of the factorization.
\end{proof}

Using the usual trick of rewriting a pair of factorizations of a
permutation as a single factorization of the identity, the following
corollary is immediate:

\begin{cor}
The parity of a factorization of a permutation into transpositions is
independent of the factorization chosen.
\end{cor}

Returning our attention to the noncrossing partition lattice, we find
that $NC_n$ is closely connected with the factorizations of an
$n$-cycle into transpositions.  First, we introduce some definitions.

\medskip
\noindent
{\bf Minimal factorizations.}  A factorization of a permutation into
transpositions is called \emph{minimal} if it has minimal length
(i.e., the smallest number of factors) among all such factorizations.
We can then define an ordering on the set of permutations by declaring
that $\sigma < \tau$ if there is a minimal factorization of $\tau$
with a ``prefix'' that is a minimal factorization of $\sigma$.  In
other words, there should be a minimal factorization $t_1 \circ t_2
\circ \cdots \circ t_k = \tau$ such that $t_1 \circ t_2 \circ \cdots
t_\ell = \sigma$ with $\ell \leq k$.

\medskip

One quick note about multiplication.  There are two natural
conventions for multiplying permutations: functional notation and
algebraist notation.  We use the algebraic convention throughout so
that $(1,2)(1,3)$ is $(1,2,3)$ rather than $(1,3,2)$.  Because of all
the symmetries of the objects under consideration this is only of
minor importance, but it means that we need to write ``prefix'' rather
than ``suffix'' in the foregoing definition.

\begin{lem}\label{lem:sym}
The poset of permutations less than or equal to the $n$-cycle
$(1,2,3,\ldots,n)$ is isomorphic with the noncrossing partition
lattice $NC_n$.
\end{lem}

\begin{proof}[Sketch of proof]
Instead of giving a complete proof we write down the isomorphism and
omit the details.  The interested reader can find a complete proof in
\cite{Br01}.  Given a noncrossing partition $\sigma$, we convert it
into a permutation $\pi(\sigma)$ by writing each block as a disjoint
cycle.  The order in which the elements occur is the natural linear
order, or stated more geometrically, we read the labels of each convex
hull clockwise.  To illustrate, the noncrossing partition $\sigma$
illustrated in Figure~\ref{fig:ncross} corresponds to the permutation
$\pi(\sigma) = (1,4,5)(2,3)(6,8)$.  Showing that the permutations
corresponding to noncrossing partitions are prefixes of reduced
factorizations of the $n$-cycle $(1,2,\ldots,n)$ is relatively easy,
as is the correspondence between the two orderings.  The only slightly
tricky part is showing that this map is onto.
\end{proof}

\begin{figure}[ht]
\includegraphics{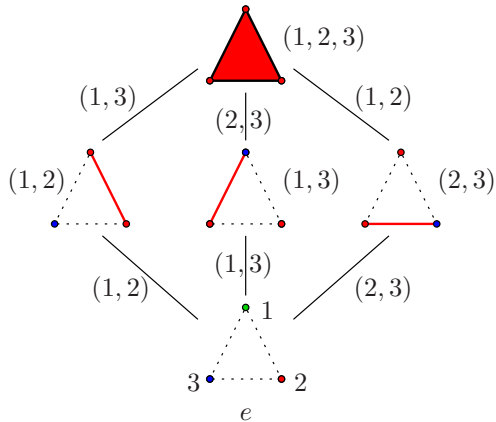}
\caption{The noncrossing partition lattice for $n=3$ with edge and
  vertex labels.\label{fig:nc3}}
\end{figure}

What Lemma~\ref{lem:sym} proves, in essence, is that the graph whose
vertices are elements of $NC_n$ and whose edges are its covering
relations corresponds to a portion of the right Cayley graph of $S_n$
with respect to the set of all transpositions.  We remind readers that
the \emph{Cayley graph} of a group $G$ with respect to a generating
set $A$ is a directed graph with vertices labeled by the elements of
$G$ and edges indexed by the set $G\times A$, where the edge $(g,a)$
connects $g$ to $g\cdot a$ \cite{GrMa64}.  Because of this
identification, not only does every element of $NC_n$ have a
permutation assigned to it, but if $\tau$ covers $\sigma$ then this
edge is labeled by a transposition (equal to
$\pi(\sigma)^{-1}\pi(\tau)$).  We have illustrated this labeling for
$NC_3$ in Figure~\ref{fig:nc3}.  Notice that the Cayley graph
interpretation also ensures that the sequence of edge labels in a
maximal chain, read from the bottom to the top, multiply together to
give $(1,2,\ldots,n)$ and that these correspond exactly to its minimal
factorizations.  For example, the minimal factorizations of $(1,2,3)$
into transpositions (in algebraist notation) are $(1,2)(1,3) =
(1,3)(2,3) = (2,3)(1,2) = (1,2,3)$ and these are the labels on the
three maximal chains seen in Figure~\ref{fig:nc3}.

More generally, multiplying together the sequence of transpositions
labeling the edges in a saturated chain connecting $\sigma$ to $\tau$
yields the permutation $\pi(\sigma)^{-1}\pi(\tau)$.  Using
Lemma~\ref{lem:sym} it is now easy to establish the following:

\begin{prop}\label{prop:words}
If $\rho < \sigma< \tau$ in $NC_n$ and
$(\alpha_1,\alpha_2,\ldots,\alpha_k)$ is a sequence of transpositions
labeling a saturated chain from $\rho$ to $\sigma$, then there is a
unique element $\sigma'$ in $NC_n$ and a saturated chain from
$\sigma'$ to $\tau$ with the exact same sequence of labels.
Similarly, if $(\beta_1,\beta_2,\ldots,\beta_k)$ is a sequence of
transpositions labeling a saturated chain from $\sigma$ to $\tau$,
then there is a unique element $\sigma''$ in $NC_n$ and a saturated
chain from $\rho$ to $\sigma''$ with this same sequence of labels.
\end{prop}

\begin{proof}
This statement is actually a consequence of the observation that the
set of transpositions in $S_n$ is closed under conjugation.  Thus, if
$\alpha$ and $\beta$ are single transpositions and $\alpha \beta$ is a
minimal factorization, then there is another minimal factorization
$\gamma \alpha$, where $\gamma$ is the transposition
$\alpha\beta\alpha = (\alpha \beta \alpha^{-1})$.  Iterating this idea
allows us to move labeled chains up and down as much as we want.
Finally, the partitions $\sigma'$ and $\sigma''$ must be unique since
we know exactly the permutations to which they correspond under the
identification with the Cayley graph.
\end{proof}

The locally self-dual property is now almost immediate:

\begin{proof}[Proof that $NC_n$ is (locally) self-dual]
Given $\sigma$ and $\tau$ in $NC_n$ such that $\sigma < \tau$, let $P$
be the poset of noncrossing partitions between them.  Define a map
$f:P\to P$ so that for each $\rho$ in $P$ the labels on a saturated
chain from $\sigma$ to $\rho$ are also the labels on a saturated chain
from $f(\rho)$ to $\tau$.  By Proposition~\ref{prop:words} a unique
such element exists, and it is easy to check that $f$ is a
well-defined order reversing isomorphism from $P$ and $P$.
\end{proof}

Finally, we note that the symmetric groups are examples of a broader
class of groups called \emph{finite reflection groups} or \emph{finite
  Coxeter groups}.  In the same way that the noncrossing partition
lattice is closely connected with the symmetric group, there is an
entire series of lattices, one for each finite Coxeter group
\cite{AtRe}, \cite{Bessis}, \cite{BrWa02}.  Each of these general
noncrossing partition lattices is graded, bounded, self-dual, and
locally self-dual, and the general proofs are essentially the same as
the ones given here.  Because of these elegant patterns, the number of
elements in these additional lattices have come to be called
\emph{generalized Catalan numbers}.  In the course of the article we
occasionally comment on properties that extend to these additional
situations.

\section{Braid groups.}
We are now ready to establish our first connection: noncrossing
partitions and the braid groups.  The braid groups are related to many
areas of mathematics including mathematical physics, quantum groups,
von Neumann algebras, and, not too surprisingly, three-manifold
topology (see, for example \cite{Bi01}, \cite{JoSu97},
\cite{KaRoTu97}, \cite{MuKu99}, \cite{PrSo97}, or \cite{Ro03}).  The
surprise is that they are also intimately related to noncrossing
partitions.

\begin{figure}[ht]
\begin{tabular}{c}
\includegraphics[scale=1]{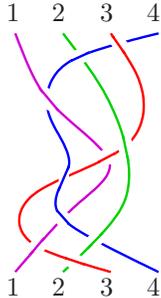}
\end{tabular}
\caption{An element of $B_4$.\label{fig:b4}}
\end{figure}

Roughly stated, a \emph{braid on $n$ strings} keeps track of how $n$
strings can be twisted in space so that each strand is a smooth
embedded monotonically decreasing curve in $\mathbb{R}^3$ (i.e., its
partial with respect to $z$ is always strictly negative).  Various
conventions need to be established, such as that strands must start
and end in some standardized configuration, that strands cannot
intersect, and that perturbations of a braid that maintain these
conventions are considered to be the same element.  The collection of
all braids is turned into the braid group $B_n$ once multiplication is
defined by attaching one braid to the top of another (see
Figure~\ref{fig:b4} for a typical element of $B_4$).  Observe that
there is a group homomorphism from $B_n$ onto $S_n$ that simply
forgets which way crossings took place and merely records the
permutation of the strings (compare, for example, Figures~\ref{fig:b4}
and~\ref{fig:sym}).  The way the noncrossing partition lattice enters
the picture is through the structure of its order complex.

Each partially ordered set $P$ can be turned into a simplicial complex
in a very simple fashion.  The vertices of the complex are labeled by
the elements of the poset, and we add a simplex corresponding to a set
of vertices if and only if the elements that label them form a chain
in $P$.  The result is called the \emph{order complex} of $P$ (or its
\emph{geometric realization}).

Let $\Delta(NC_n)$ denote the order complex of $NC_n$.  Notice that
the one-cells in $\Delta(NC_n)$ correspond to two-element chains in
$NC_n$.  In other words, they correspond to pairs of noncrossing
partitions $\sigma$ and $\tau$ with $\sigma<\tau$.  We can carry over
the labeling from $NC_n$ to the one-cells in $\Delta(NC_n)$ as
follows: we label the oriented edge from $\sigma$ to $\tau$ in
$\Delta(NC_n)$ by the permutation $\pi(\sigma)^{-1}\pi(\tau)$.  Our
final step is identify certain simplices in $\Delta(NC_n)$ with each
other.  The rule is that two simplices are identified if and only if
they have identically labeled, oriented one-skeletons (and they are
identified so that these labels match, of course).  The resulting
quotient is a complex we call the \emph{Brady-Krammer complex} $BK_n$,
since it was discovered independently by Tom Brady and Daan Krammer.
What Brady and Krammer proved in \cite{Br01} and \cite{Kr00},
respectively, was that this procedure results in a complex whose
fundamental group is the braid group $B_n$ and whose universal cover
is contractible.  In other words, they proved the following result:

\begin{thm}[{\bf Brady, Krammer}]\label{thm:bk}
The complex $BK_n$ is an Eilenberg-Maclane space for the braid group
$B_n$.
\end{thm}

As in the previous section, the procedure described extends naturally
to the general noncrossing partition lattices associated with the
other finite reflection groups.  In each case the resulting complex is
an Eilenberg-Maclane space for a group, and this group is related to
the finite reflection group in the same way that the braid group is
related to the symmetric group.  These other groups are called
\emph{finite-type Artin groups}.  See \cite{Bessis} or \cite{BrWa02}
for a proof of this extension or \cite{Be99}, \cite{Ch95},
\cite{ChDa95}, or \cite{CoWa02} for more about finite-type Artin
groups.

\section{Parking functions.}
\begin{figure}[th]
\includegraphics[width=4in]{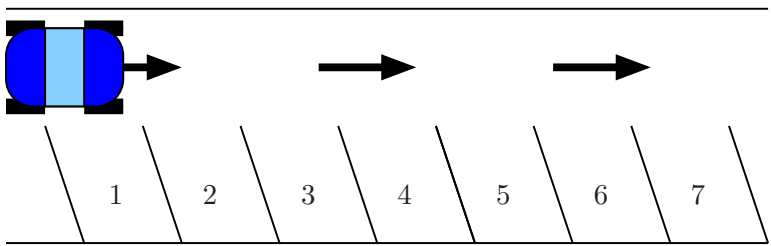}
\caption{A one-way street with $7$ parking spots.\label{fig:park}}
\end{figure}

For our second illustration we shift to a classic problem from
combinatorics.  Imagine a sequence of $n$ cars entering a one-way
street one at a time with $n$ parking spots available, as shown in
Figure~\ref{fig:park}.  Each driver has a preferred parking spot and
attempts to park there first.  Failing that, he or she parks in the
next available space.  If any of the drivers is forced out of the
street, then this sequence of preferences has failed.  We refer to
any sequence of preferences that enables all $n$ of the drivers to
park successfully as a \emph{parking function}.

As an initial observation, it is easy to see that if two or more
drivers prefer the final parking space, the sequence fails.
Similarly, if three or more drivers prefer either of the last two
spaces, the sequence fails, etc.  Perhaps surprisingly, these are the
only restrictions.  Thus, an equivalent definition of a parking
function would describe it as a sequence $(a_1,\ldots,a_n)$ of
positive integers whose rearrangement (and relabeling) as a
nondecreasing sequence $b_1\leq b_2 \leq \ldots \leq b_n$ satisfies
the inequalities $b_i\leq i$ for each $i$ in $[n]$.  From this altered
definition it is not hard to show that the number of parking functions
is $(n+1)^{n-1}$.  Combinatorialists, of course, recognize this as the
number that counts labeled rooted trees on $[n]$ (or, equivalently,
the number of acyclic functions on $[n]$).  In order to see the
connection with the noncrossing partition lattice, consider the
following definition due to Richard Stanley.

Suppose that $\tau$ covers $\sigma$ in $NC_{n+1}$, that $B$ and $B'$
are the two blocks of $\sigma$ that combine to form a block in $\tau$,
and without loss of generality, that $\min B <\min B'$.  We define a
new label on this covering relation by the largest element of $B$ that
is below each element of $B'$.  Using this edge-labeling Stanley was
able to show the following \cite{St97}:

\begin{thm}[{\bf Stanley}]\label{thm:parking}
The labels on the maximal chains in $NC_{n+1}$ are exactly the parking
functions of length $n$, each occurring once.
\end{thm}

\noindent
The reader wishing to read more about parking functions and
noncrossing partitions might want to consult \cite{Bi02}, \cite{Ka02},
\cite{StPi02}, \cite{Ya00}, or especially, \cite{St97}.

\section{Free probability.}
Our third sighting of the noncrossing partition lattice is in a
noncommutative version of probability. Because of the nature of the
subject matter, the discussion in this section is less detailed than
in previous ones.  Readers wishing to read more about the connection
between free probability and noncrossing partitions should probably
begin with the excellent survey article by Roland Speicher
\cite{Sp97}.  We start with a brief discussion of classical
probability.

Let $X$ be a random variable having a probability density function
$f(x)$.  For the reader unfamiliar with probability theory, a good
example of a probability density function is a nonnegative continuous
function from $\mathbf{R}$ to $\mathbf{R}$ whose integral over the
reals is $1$.  The \emph{expectation} of a function $u(X)$ is then
defined as
\[E(u(X)) = \int_{-\infty}^\infty u(x) f(x) dx.\]

The first expectations students usually encounter are the mean $\mu =
E(X)$ and the variance $\sigma^2 = E((X-\mu)^2)$.  In a mathematical
probabilty and statistics course they meet the higher moments $E(X^n)$
as well as the moment generating function
\[M(t) = E(e^{tX}) = \sum_{n\geq 0} E(X^n) \frac{t^n}{n!},\]
\noindent
which allows the calculation of all of the moments by evaluating a
single integral.  The coefficents in the moment generating function
are called the (classical) \emph{moments} of $X$.  The coefficients of
$\log M(t)$ are called the (classical) \emph{cumulants} of $X$.  The
main advantage of the cumulants is that they contain the same
information as the moments but that the $n$th cumulant of the sum of
two random variables $X$ and $Y$ is the sum of their $n$th
cumulants--provided they are independent.

Before launching into the ``noncommutative'' version, we should say a
brief word about noncommutative geometry in general.  Noncommuatative
geometry is a philosophy whereby standard geometric arguments on
topological spaces are converted into algebraic arguments on their
commutative $C^*$-algebras of functions.  The motivation comes from
mathematical physics and the need to integrate quantum mechanics
(which is noncommutative) with classical physics.  The main
observation is that there is a nice correspondence between
``reasonable'' topological spaces $X$ and the collections of
continuous maps from these spaces to the complex numbers.  Each such
collection has a structure known as a $C^*$-algebra.  In fact, the
correspondence is strong enough that a space $X$ can be recovered from
the commutative $C^*$-algebra to which it gives rise.  The philosophy,
in short, is to reformulate each concept from classical topology,
geometry, calculus of manifolds, and so forth in terms of properties
of $C^*$-algebras and then to turn these equivalent formulations into
definitions for a ``noncommutative'' version of this concept.  Even
though there are no longer any topological spaces or points or open
sets---only $C^*$-algebras---the classical structures can be used to
develop intuition and provide a guide to the types of theorems and
results that should be expected from the noncommutative world.  To
date this program has been remarkably successful.  The original book
by Alain Connes \cite{Co94} or the more recent (and shorter) survey
articles such as \cite{Co00} are a good place for an inexperienced
reader to begin.

Returning to probability theory, we remark that researchers have
defined a noncommutative probability space to be a pair
$(\mathcal{A},\phi)$, where $\mathcal{A}$ is a complex-unital-algebra
equipped with a unital linear functional $\phi$ called expectation.
There is also a noncommutative version of independence known as
``freeness.''  Without getting into the details, the combinatorics of
noncrossing partitions is very closely involved in the noncommutative
version of cumulants.  In fact, some researchers who study free
probability have described the passage from the commutative to the
noncommutative setting of probability as a transition from the
combinatorics of the partition lattice $\Pi_n$ to the combinatorics of
noncrossing partitions $NC_n$.  The article by Roland Speicher
\cite{Sp97} is an excellent exposition of this topic.  To give just
one hint at the underlying argument, we note that new counting
problems that involve summing up an old counting problem over all of
the possible partitions of that problem into smaller problems of the
same type often lead to solutions that contain exponentials (see
\cite[chap. 5]{EC2} for a precise development of this theme).  Thus
the exponential function in the integral defining the moment
generating function gets converted into a sum over the elements in the
partition lattice.  In the noncommutative context, the crossing
partitions are prevented from playing a role, so the sum takes place
over the noncrossing partition lattice instead.

\section{Summary.}
As we have seen, the lattice of noncrossing partitions might surface
in any situation that involves (1) the symmetric groups, (2) the braid
groups, (3) free probability, or (4) the Catalan numbers (say, in
conjuction with the combinatorics of trees or the combinatorics of
parking functions).  They also show up in real hyperplane
arrangements, Pr\"ufer codes, quasisymmetric functions, and Hopf
algebras, but detailing all of these connections would lead us too far
afield.  (see \cite{Ai01}, \cite{At}, \cite{At98}, \cite{AtLi99},
\cite{AvBe03}, \cite{Eh96}, \cite{EC1}, \cite{EC2}, or \cite{St97} for
details.)  We conclude with one final illustration: the associahedron.

\begin{figure}[ht]
\begin{tabular}{ccc}
\begin{tabular}{c}
\includegraphics{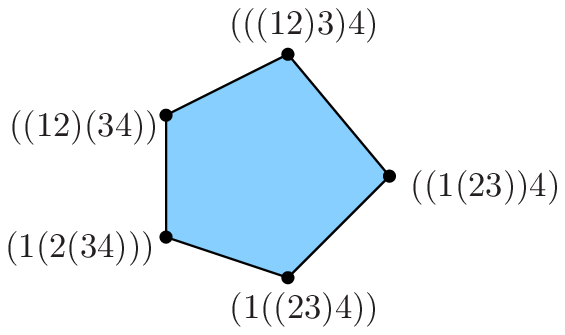}
\end{tabular} & \hspace*{1cm} & \begin{tabular}{c}
\includegraphics[scale=.5]{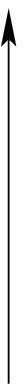}
\end{tabular}\end{tabular}
\caption{The two dimensional associahedron.\label{fig:assoc2}}
\end{figure}

\begin{figure}[ht]
\begin{tabular}{cc}
\begin{tabular}{c}
\includegraphics[width=2in]{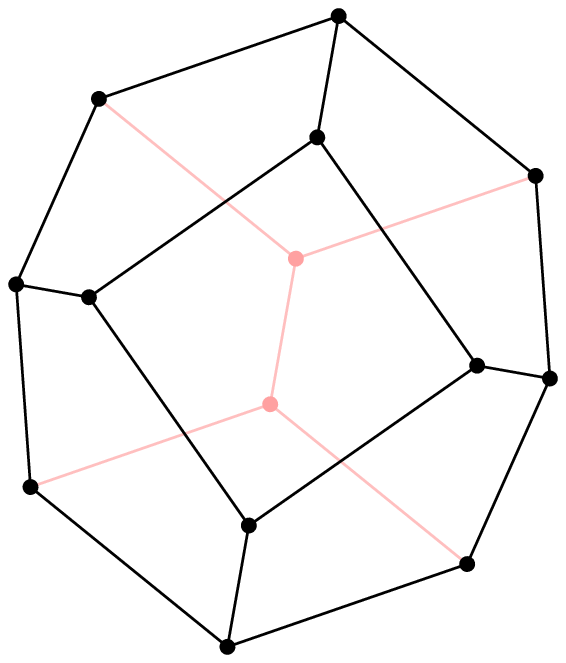}\end{tabular} & \begin{tabular}{c}
\includegraphics{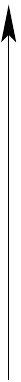}\end{tabular}
\end{tabular}
\caption{The $3$-dimensional associahedron.\label{fig:assoc3}}
\end{figure}

\medskip
\noindent
{\bf Associahedron.}  As we remarked earlier, the Catalan numbers
count the number of ways to associate a list of numbers.  If we also
consider partial associations such as $(1 2 3)(4 5)$, we get a partial
ordering of partial associations (where removing a pair of matched
parentheses corresponds to moving up in the ordering).  This partial
ordering has been shown to be the lattice of faces for a convex
polytope known as the \emph{associahedron} (or \emph{Stasheff
polytope}).  For example, the partial associations of four elements
form a pentagon (Figure~\ref{fig:assoc2}).  The full associations
label the vertices as indicated, and the partial associations label
the edges.  For example, $((1 2) 3 4)$ is the partial association
corresponds to the edge connecting the vertices labeled $(((1 2) 3)
4)$ and $((1 2)(3 4))$ since it can be obtained from either one by
removing a pair of parentheses.  Similarly, the fourteen associations
listed in Figure~\ref{fig:associate} can be identified with the
fourteen vertices of the polytope shown in Figure~\ref{fig:assoc3}.

It is the ``Morse theory'' of this polytope that has a close
connection with the noncrossing partition lattice.  If we chose a
height function (i.e., a linear map from the Euclidean space
containing the polytope onto the reals) so that none of the edges are
horizontal (i.e., none of the direction vectors of the edges lie in
the kernel of this map), then at each vertex we can count the number
of edges pointing ``up'' and the number pointing ``down.''  For the
pentagon pictured in Figure~\ref{fig:assoc2} (with a height function
that orthogonally projects onto the arrow shown) there is one vertex
with both edges pointing up, three vertices with one up and one down,
and one vertex with both edges pointing down.  We can tally these
results in a vector known as the \emph{$h$-vector} of the polytope.
Thus, the $h$-vector of the pentagon is $(1,3,1)$.  Although it is not
obvious from this definition, this sequence is independent of the
height function chosen.

The rank function for the noncrossing partition lattice $NC_3$ (i.e.,
the number of elements it contains at each height) can also be
summarized by a vector.  As seen in Figure~\ref{fig:nc3}, the vector
for $NC_3$ is $(1,3,1)$ since there is one element with height zero,
three elements with height one, and one element with height two.  The
recurrence of the vector $(1,3,1)$ is not a coincidence.  Notice that
the associahedron shown in Figure~\ref{fig:assoc3} has an $h$-vector
$(1,6,6,1)$ that corresponds nicely with the rank function for $NC_4$
(Figure~\ref{fig:nc4}).
 
Fomin and Zelevinsky have recently defined for each finite
crystallographic reflection group a Euclidean polytope known as a
\emph{generalized associahedron} \cite{FoZe}, \cite{ChFoZe02}, and in
each case the $h$-vector for the general polytope matches the rank
function for the corresponding lattice \cite{At}.  These lattices and
polytopes, and the observed connections between them, are only a few
years old at this point. They are the subject of many ongoing research
projects.  I am sure that we will be discovering additional remarkable
properties of these objects for many years to come.

\vspace{1em}
\noindent
\textbf{ACKNOWLEDGMENTS.} In 1999 Texas A\&M University hired three
new assistant professors: myself (geometric group theory), Catherine
Huafei Yan (combinatorics), and Ken Dykema (free probability).  The
discovery that we were all interested in the structure of the
noncrossing partition lattice from vastly different perspectives
marked the early beginnings of this article.  The more proximate cause
was a general-interest talk I gave at the conference celebrating Jim
Cannon's sixtieth birthday held in Park City, Utah, in June, 2003.  I
would like to thank the organizers for their invitation.  I would also
like to thank Christos Athanasiadis, Tom Brady, Ken Dykema, John
Meier, Vic Reiner and Jim Stasheff for comments on early versions of
the article.  As a final note, a colorful version of this article is
available from my homepage ({\tt
http://www.math.ucsb.edu/$\sim$mccammon/}).  The author was partially
supported by the National Science Foundation.

\def\cprime{$'$}

\medskip
\noindent
\textbf{JON MCCAMMOND} is an associate professor at the University of
California, Santa Barbara.  Although his research is primarily focused
on geometric group theory, he has a more than passing interest in
geometric combinatorics (particularly polytopes), theoretical computer
science (particularly finite state automata), and differential
geometry (particularly symmetric spaces).

\end{document}